\documentclass[preprint,12pt]{elsarticle}
\usepackage{lipsum}
\usepackage{graphicx}
\usepackage{amsmath,leftidx,amsthm}
\usepackage[all]{xy}
\usepackage{xspace}
\usepackage{amssymb}
\usepackage[latin1]{inputenc}
\usepackage{graphicx,color}
\usepackage{hyperref,fancyhdr}
\usepackage{fullpage}
\usepackage{versions}
\usepackage{textcomp,listings}

\newtheorem{theorem}{Theorem}[section]
\newtheorem{lemma}[theorem]{Lemma}
\newtheorem{cor}[theorem]{Corollary}
\newtheorem{conj}[theorem]{Conjecture}
\theoremstyle{definition}
\newtheorem{definition}[theorem]{Definition}

\newtheorem{proposition}[theorem]{Proposition}

\theoremstyle{remark}

\numberwithin{equation}{section}

  \renewcommand{\Re}{\mbox{\rm Re}}
  \renewcommand{\Im}{\mbox{\rm Im}}
    
   \newcommand{\ph}{\mbox{$\varphi$}}
    
    \renewcommand{\phi}{\varphi}
   \newcommand{\Ht}{\mbox{$H^{2}$}}
   
   \newcommand{\Hi}{\mbox{$H^{\infty}$}}
   \newcommand{\D}{\mbox{$\mathbb{D}$}}
   \newcommand{\C}{\mbox{$C_{\varphi}$}}
   
   \newcommand{\T}{\mbox{$T_\psi$}}
   
   \newcommand{\W}{\mbox{$W_{\psi,\varphi}$}}
   \newcommand{\Ws}{\mbox{$W_{\psi,\varphi}^{\textstyle\ast}$}}
 	\newfont{\caps}{cmcsc9}  
 	\newfont{\jour}{cmti9}  
\theoremstyle{remark}

\definecolor{green}{rgb}{0,0.5,0}
\definecolor{dkgreen}{rgb}{0,0.6,0}
\definecolor{gray}{rgb}{0.5,0.5,0.5}
\definecolor{mauve}{rgb}{0.58,0,0.82}
\definecolor{purple}{rgb}{0.58,0,0.62}

\renewcommand{\phi}{\varphi}
\renewcommand{\Re}{\mbox{Re\,}}
\renewcommand{\Im}{\mbox{Im\,}}

\def\and{{\quad\text{and}\quad}}

\journal{Journal of Mathematical Analysis and Applications }

\begin{document}

\begin{frontmatter}



\title{Spectra of Weighted Composition Operators with Quadratic Symbols}
\tnotetext[t1]{The authors were funded by an NSF-CURM grant.}

\author[n]{Jessica Doctor} 
\author[n]{Timothy Hodges}
\author[t]{Scott Kaschner}
\author[n]{Alexander McFarland}
\author[n]{Derek Thompson\corref{cor1}}
\address[t]{Department of Mathematics, 
Butler University, 4600 Sunset Ave, Indianapolis, IN 46208}
\address[n]{Department of Mathematics, Taylor University, 236 W. Reade Ave, Upland, IN 46989}
\ead{theycallmedt@gmail.com}

\begin{abstract}
   Previously, spectra of certain weighted composition operators \W\ on \Ht\ were discovered under one of two hypotheses: either \ph\ converges under iteration to the Denjoy-Wolff point on all of \D\ rather than compact subsets, or \ph\ is ``essentially linear fractional''. We show that if \ph\ is a quadratic self-map of \D\ of parabolic type, then the spectrum of \W\ can be found when these maps exhibit both of the aforementioned properties, as most of them do. 
\end{abstract}

\end{frontmatter}

\section{Introduction} 

Let \Ht\ denote the classical \textit{Hardy space}, the Hilbert space of analytic functions $\displaystyle f(z) = \sum_{n=0}^{\infty}a_n z^n$ on the open unit disk \D\ such that $$\|f\|^{2}=\sum_{n=0}^{\infty}|a_n|^{2}<\infty.$$ A \textit{composition operator} \C\ on \Ht\ is given by $\C f = f \circ \ph$. We call \ph\ the \textit{symbol} of the associated composition operator. When \ph\ is an analytic self-map of \D, the operator \C\ is bounded. Composition operators on \Ht\ have been extensively studied for several decades; see \cite{cowen1995composition} and \cite{shapirocomposition} for seminal books on the subject. One reason for their study is the deep connection to \textit{Toeplitz operators}: $T_\psi$ on $H^2$ is given by $T_\psi f = P(\psi f)$ where $P$ is the projection back to $H^2$. When $\psi \in \Hi$, the space of bounded analytic functions on \D, we simply have $T_\psi f = \psi f$, since $\psi f$ is guaranteed to be in \Ht, and all such Toeplitz operators are bounded. Throughout this paper, we will assume $\psi \in \Hi$.  We write $\W := T_\psi \C$ and call such an operator a \textit{weighted composition operator}. 

Another reason for the study of composition operators is the interesting interplay between the operator-theoretic behavior of \C\ and the behavior of \ph\ on \D. The Denjoy-Wolff theorem guarantees that any analytic self-map of \D\  (apart from elliptic automorphisms) will have a unique attracting fixed point $w$ in $\overline{\D}$, and \ph\ converges, uniformly on compact subsets of \D, to that point under iteration. This behavior generally splits into three categories, each with different results for the behavior of \C: either \ph\ has a fixed point properly in \D\ with $|\ph'(w)| < 1$ (elliptic),  or $w$ is on the boundary (in a radial limit sense; \ph\ need not be analytic on $\partial \D$), and $\ph'(w) < 1$ (hyperbolic) or $\ph'(1) = 1$ (parabolic). When the fixed point is on $\partial \D$, we define $\ph'(w)$ as an angular derivative, which coincides with the typical derivative when \ph\ is analytic on $\partial \D$. We will only be considering maps for which \ph\ is analytic on all of $\overline{D}$. Generally, the parabolic case is considered the most difficult  of the three (see, e.g., \cite[Chapter 7]{cowen1995composition}.)

The \textit{spectrum} of an operator $T$ on a Hilbert space $H$, denoted $\sigma(T)$, is given by $\{ \lambda \in \mathbb{C} : T - \lambda I \textrm{ is not invertible} \}$. The essential spectrum  $\sigma_{e}(T)$ is the spectrum of $T$ in the Calkin algebra $\mathcal{B}(H)/ \mathcal{K}(H)$, the bounded operators on a Hilbert space $H$ modulo the compact operators. The essential spectrum is always a subset of the spectrum. Spectra of weighted composition operators on \Ht\ were studied in \cite{cowen2016spectra} and \cite{bourdon2012spectra}; both papers depended on the behavior of \ph\ on \D. In \cite{cowen2016spectra} it was assumed that \D\ converges uniformly to the Denjoy-Wolff point on \textit{all} of \D\ rather than compact subsets (denoted UCI for \textit{uniformly convergent iteration}). In \cite{bourdon2012spectra}, \ph\ was assumed to be \textit{essentially linear fractional}, a generalization of linear fractional maps, which play nicely with the various standard tools used to study composition operators on \Ht\ (e.g. the reproducing kernels on \Ht\ are linear fractional.)  Both papers were able to completely characterize $\sigma(\W)$ if $\psi \in \Hi$ and \ph\ is hyperbolic or elliptic, and both give partial results if \ph\ is parabolic. 

Therefore, we are interested in maps for which \ph\ is parabolic, so that we may combine the results of these prior publications. In this paper, we study the spectrum of \W\ when \ph\ is a quadratic self-map of \D\ of parabolic type.  For most of these maps, we give a complete characterization of $\sigma(\W)$. 

In the next section, we establish preliminary facts about \ph\ in this setting and divide it into two cases. In Section 3, we establish that \ph\ exhibits UCI in both cases. In Section 4, we characterize when \ph\ is essentially linear fractional and then give our main spectral results in Section 5. We end with some commentary and questions about further work in this direction in Section 6.

\section{\ph\ is of parabolic type}

Throughout this paper we assume that \ph\ is analytic in $\overline{D}$ with a single fixed point $w$ on the boundary, with $\ph'(w) = 1$. Since $C_{e^{i\theta }z}$ is a unitary operator, we can conjugate \C\ by $C_{e^{i\theta }z}$  to rotate the fixed point without affecting the spectrum; therefore we may assume that $w = 1$. We can then simplify the generic quadratic $\ph(z) = a_2 z^2+a_1 z+a$ by noting 

\begin{align}
    a_2+a_1+a= 1 \\
    2a_2+a_1 = 1
\end{align}

which gives $a_1 = 1 - 2a_2$, and $a_2 = a$, so we have $\ph(z) = az^2+(1-2a)z+a$.

We now need to confirm exactly when \ph\ maps \D\ into \D. 

\begin{proposition}\label{dtod}
Suppose $\ph(z) = az^2+(1-2a)z+a$ so that $\ph(1) =1, \ph'(1) = 1$. Then \ph\ is a self-map of \D\ if and only if $|a-\frac{1}{4}| \leq \frac{1}{4}$. Furthermore, the only elements of $\partial \D$ that are not mapped into $\D$ are $\ph(1)=1$, and when $|a-\frac{1}{4}| = \frac{1}{4}$,  $\ph(-1) = 4a-1$. 
\end{proposition}

\begin{proof}
Note that $\ph(-1) = 4a - 1$, so we require that $|4a-1| \leq 1$. Equivalently, we have $|a-\frac{1}{4}| \leq \frac{1}{4}$, so that $a$ is contained in the circle of radius $\frac{1}{4}$ centered at $\frac{1}{4}$.

Now suppose instead that $|a-\frac{1}{4}| \leq \frac{1}{4}$ is given. Then, we will find the image of the unit circle under $\ph$. Thinking of $\ph$ as a function of $a$, we know by the Maximum Modulus Principle that $|f(z)|$ will be maximized when $|a-\frac{1}{4}|=\frac{1}{4}$. Therefore, writing $a = A+Bi$, we have $A^2+B^2 = \frac{1}{2}A$ and if $z = X+Yi$, $X^2+Y^2 = 1$. 

If $\ph(z) = az^2+(1-2a)z+a = (A+Bi)(X+Yi)^2+(1-2A-2Bi)(X+Yi)+(A+Bi)$, then $|\ph(z)|^2$ can be directly computed, resulting in a 28-term algebraic expression we will not suffer the reader to witness. However, using the aforementioned equations, we are able simplify the expression to $   2AX^2 - 2 A    + 1 $. This is a real-valued quadratic in $X$ defined on $[-1,1]$, with vertex $(0,-2A+1)$ (which has a nonnegative $y-$value for $ 0 \leq A \leq \frac{1}{2}$, as needed). Since the graph is a parabola whose defining equation has a positive leading coefficient, the maximum values in the domain are found at the endpoints $x = 1$ and $x = -1$, which both result in $|\ph(z)|^2 = 1$. Note, however, equality is always attained at $\ph(1)=1$, but equality at $x=-1$ only happens for $\ph(-1)=4a-1$ when $|a-\frac{1}{4}|=\frac{1}{4}$; otherwise $\ph(-1)$ is properly in \D. In any case, for any other value of $X$, we have $|\ph(z)|^2 < 1$, so \D\ is mapped by \ph\ into \D\ as desired. 
\end{proof}

While the primary purpose of Proposition \ref{dtod} is discovering exactly when \ph\ is a self-map of \D, we should also note a few key facts from this result that we will use later. In particular, we see that $\Re \textrm{ } a > 0$ and $|a| \leq \frac{1}{2}$. Furthermore, we will take advantage of knowing exactly which points on the unit circle are mapped by \ph\ back onto the unit circle. As we have seen in this proof, $|a-\frac{1}{4}| < \frac{1}{4}$ is a separate case from $|a-\frac{1}{4}| = \frac{1}{4}$. This will remain true throughout the paper.

\section{\ph\ converges uniformly under iteration on all of \D\ (UCI)}

Rather than iterating our quadratic maps directly, we will circumvent this issue by taking alternate approaches to showing that $\ph_n \rightarrow 1$ uniformly on all of \D. For $|a-\frac{1}{4}| < \frac{1}{4}$, we rely on  arguments from complex dynamics. When $|a-\frac{1}{4}| = \frac{1}{4}$, we work directly with the following  alternate characterization of uniform convergence.

\begin{proposition}
A sequence of functions $f_n$ converges uniformly to its pointwise limit $f$ on all of a domain $D$ if and only if $$\lim_{n \rightarrow \infty} \sup_{z \in \D} |f_n(z)-f(z)| \rightarrow 0. $$
\end{proposition}

\begin{theorem}\label{uci1}
Suppose $\ph(z) = az^2+(1-2a)z+a$ and $|a-\frac{1}{4}| = \frac{1}{4}$, $a \neq 0$. Then the iterates of \ph\ converge uniformly to $1$ on the entire open disk \D. 
\end{theorem}
\begin{proof}
We will find a recursive pattern for 
$\sup_{z \in D} |\ph_n(z)-1|$ that approaches $0$.

Suppose $\sup_{z \in \D} |\ph_{n}(z)-1| = r$ and consider $|\ph_{n+1}(z)-1|$. This factors as $$  |a||\ph_{n}(z)-1||\ph_{n}(z)-(1-\frac{1}{a})|.$$ We know that $\ph_{n}(z)$ lies on the circle with center $1$ and radius $r$ (call this circle $C$), and wish to find an upper bound on the distance of $\ph_n(z)$ from $z_0 = 1-\frac{1}{a}$. To do this, we also need to find the image of the circle $|a-\frac{1}{4}| = \frac{1}{4}$ (ignoring $a = 0$) under this transformation.  Again, if $a = A+Bi$, then $A^2+B^2 = \frac{1}{2}A$. Then we have

$$\frac{1}{a} = \frac{1}{A+Bi} = \frac{A-Bi}{A^2+B^2} = \frac{2(A-Bi)}{A} = 2-\frac{B}{A}i, $$

\noindent so we have $z_0 = 1-\frac{1}{a} = -1 + \frac{B}{A}i$, and we know $0 < A \leq \frac{1}{2}$. Without loss of generality, assume $B$ is positive. Then the furthest point from on $C$, still within $\overline{\D}$, from $z_0$ is the point where $C$ intersects the unit circle in the fourth quadrant, say $z_1$ (see Figure \ref{circle}). Writing $x^2+y^2=1$ and $(x-1)^2+y^2=r^2$, we can write $z_1$ in terms of $r$: $\left(1-r^2/2\right) -\left(r\sqrt{1-\frac{1}{4}r^2}\right)i. $ We now wish to find $|z_1 - z_0|:$

$$\left|\ph_{n}(z)-(1-\frac{1}{a})\right| \leq |z_1 - z_0| = \sqrt{1+2r^{2}+\frac{B^2}{A^2} -\frac{2B}{A} r\sqrt{1-\frac{1}{4}r^2}}. $$

To incorporate the rest of our original expression for $|\ph_{n+1}(z)-1|$, we recall that $|a| = \sqrt{A^2+B^2} = \sqrt{\frac{1}{2}A}$ and $|\ph_{n}(z)-1|=r$. We can also substitute $B = \sqrt{\frac{1}{2}A-A^2} $ to get

$$\sup_{z\in\D} |\ph_{n+1}(z)-1| \leq r \sqrt{r^{2}A+\frac{1}{4}-r\sqrt{\left(\frac{1}{2}A-A^2\right)\left(1-\frac{1}{4}r^2\right)}}.$$

\begin{figure}[h]
\includegraphics[width=0.75\textwidth]{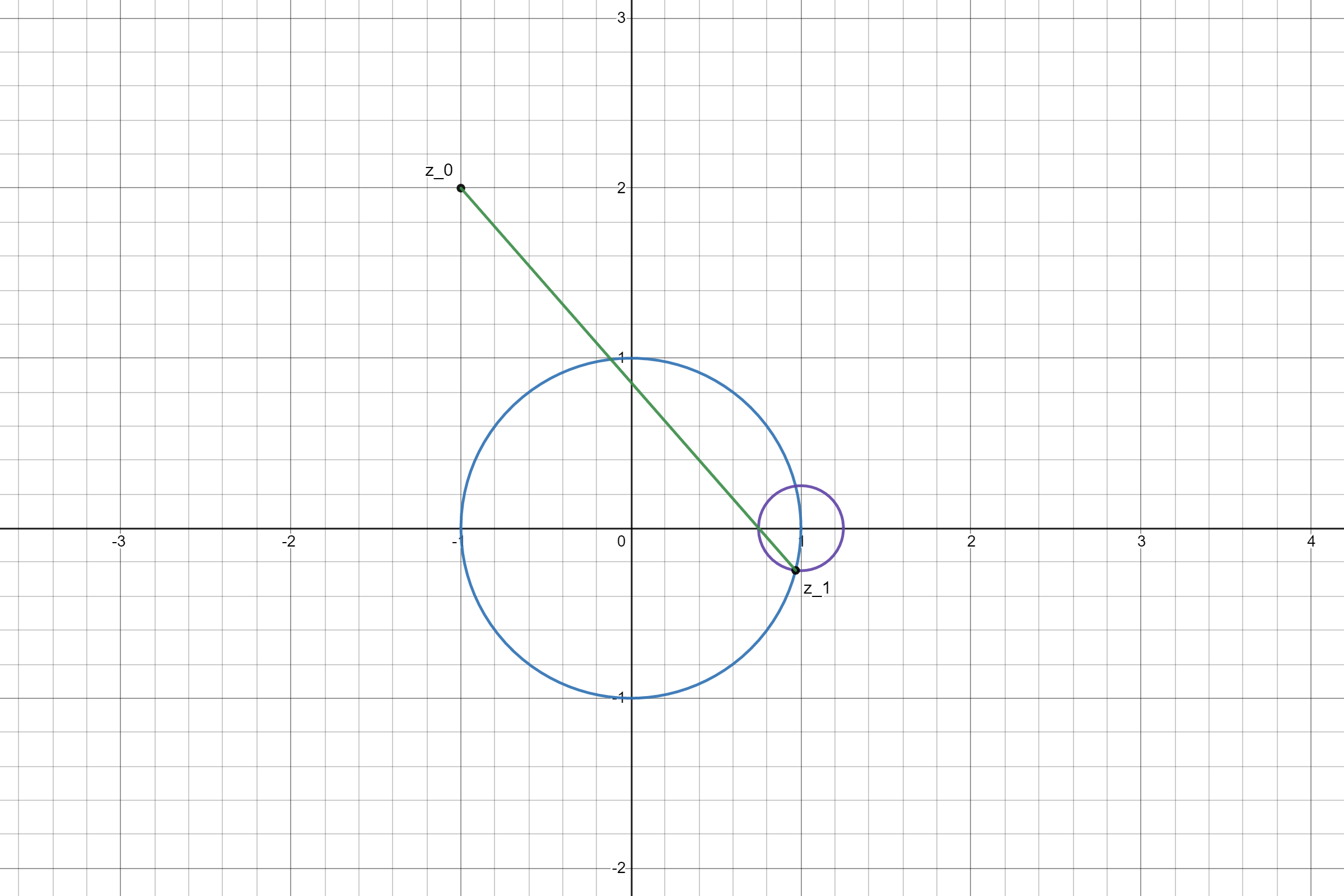}
\caption{An example scenario for $z_1$ and $z_0$.}\label{circle}
\end{figure}

The value of $r$ actually depends on which iterate of \ph\ we are considering, so we can rename it $r_n$. For our base case, consider $r_1 = |\ph(z)-1|$. Then for $z = X+Yi \in \partial \D$, if $a = A+Bi$, by direct computation we have $$|\ph(z)-1|^2 = 2-2AX^{2}+4AX-2A+4(X-1)\sqrt{A/2-A^{2}}\sqrt{1-X^{2}}-2X.$$ For smaller values of $A$, this value is maximized at $X = -1$. However, when $|(4a-1)-1|^2=|4A-2+4Bi|^2=1$, we have $A = \frac{3}{8}$, and the situation changes. For $\frac{3}{8} \leq A \leq \frac{1}{2}$, we have $ |\ph(z)-1| \leq 1 $ for all $z \in \overline{\D}$. For $0 < A < \frac{3}{8}$, the maximum value of $|\ph(z)-1|$ is $2\sqrt{1-2A}$, so $r_1 = \max \{ 2\sqrt{1-2A}, 1 \}$. 

For $0 < A \leq \frac{1}{2}$, define $f(r) = \sqrt{r^{2}A+\frac{1}{4}-r\sqrt{\left(\frac{1}{2}A-A^2\right)\left(1-\frac{1}{4}r^2\right)}}$, so $f$ satisfies $0 < f(r) < 1$ for $0 \leq r \leq \max \{ 2\sqrt{1-2A}, 1\}$. Then $r_{n+1} = r_n \sqrt{r_n^{2}A+\frac{1}{4}-r_n\sqrt{\left(\frac{1}{2}A-A^2\right)\left(1-\frac{1}{4}r_n^2\right)}} < r_n$, and thus, our sequence is decreasing. Since $r$ and $f(r)$ are both positive, we know that $r_n$ is bounded below by $0$. Since our sequence is monotonically decreasing and bounded below, it converges. Then $c = \lim r_{n+1} = \lim {r_n}$, which means that $c = cf(c)$ where $f$ is defined as above. We have already established that $f$ is positive, so the only solution is $c = 0$. Therefore,  $\lim_{n \rightarrow \infty} \sup_{z \in D} |\ph_n(z)-1| = 0$, so $\ph_n$ converges under iteration to $1$ uniformly on \textit{all} of \D.

\end{proof}

For $|a-\frac{1}{4}| < \frac{1}{4}$, we now turn to traditional results in complex dynamics, via Beardon \cite{beardon}. 

\begin{theorem}\label{uci2}
Suppose $\ph(z) = az^2+(1-2a)z+a$ and $|a-\frac{1}{4}| < \frac{1}{4}$. Then the iterates of \ph\ converge uniformly to $1$ on the entire open disk \D. 
\end{theorem}

\begin{proof}
Consider the family of degree two polynomial maps $\phi\colon\mathbb D\rightarrow\mathbb D$ given by
\[\phi(z)=az^2+(1-2a)z+a,\]
where $|a-1/4|<1/4$.  These maps all have a parabolic fixed point at $z=1$ and as maps of $\mathbb C$ are conjugate to the map $z\mapsto z-z^2$.  Specifically, for the maps $g,\sigma\colon\mathbb C\rightarrow\mathbb C$ given by \begin{eqnarray*}
g(z)&=&z-z^2\quad\mbox{and}\\
\sigma(z)&=&-\frac{1}{a}z+1,
\end{eqnarray*}
we have $g=\sigma^{-1}\circ\phi\circ\sigma$ on $\mathbb C$.  Since $\phi$ is forward invariant on the disk, $g$ is forward invariant on $\mathbb D_a:=\sigma^{-1}(\mathbb D)$, and we have the following commuting diagram
\begin{equation}\label{eq_comm_sq}\xymatrix {\relax
\mathbb D_a \ar[r]^{g} \ar[d]_{\sigma} & \mathbb D_a\ar[d]^{\sigma} \\
\mathbb D \ar[r]
    _{\phi} & \mathbb D}
\end{equation}
In particular, $\phi_n$ converges uniformly to 1 on $\mathbb D$ if $g_n$ converges uniformly to 0 on $\mathbb D_a$.

For specific details regarding the dynamics of one complex variable, we refer the reader to (\cite{beardon, milnor, carleson}).  The domain of $g$ is partitioned into two totally invariant sets, the Julia set, denoted $J(g)$, and Fatou set, denoted $\mathcal F(g)$.  The Fatou set is the set of points for which the sequence of iterates forms a normal family, and the Julia set is the complement of the Fatou set.  In this case, $\mathcal F(g)$ is just the disjoint union of the two open sets:
\begin{eqnarray*}
B(g,0)&=&\{z\in\mathbb C\colon g_n(z)\rightarrow0\}\quad\mbox{and}\\
B(g,\infty)&=&\{z\in\mathbb C\colon g_n(z)\rightarrow\infty\},
\end{eqnarray*}
called the basin of zero and the basin of infinity, respectively.  $J(g)$ has no interior, 
so it must be that $\mathbb D_a$ is a subset of $B(g,0)$ or $B(g,\infty)$.  Note that $\sigma^{-1}(0)=a$, and we have assumed that $|a-1/4|<1/4$.  Define
\[A=\{a\in\mathbb C\colon |a-1/4|<1/4\}.\]
Since
\[\left|g(a)-\frac{1}{4}\right|=\left|a-a^2-\frac{1}{4}\right|=\left|a-\frac{1}{2}\right|^2\leq\left(\left|a-\frac{1}{4}\right|+\frac{1}{4}\right)^2<\frac{1}{4},\]
we have that $A$, an open set with nonempty interior, is forward invariant by $g$, so $A$ is a subset of either $B(g,0)$ or $B(g,\infty)$.  It is easily verified that for real $a\in A$, $g_n(a)\rightarrow0$.  Thus, we have $A\subset B(g,0)$ (See Figure \ref{FIG}), so $\sigma^{-1}(0)\in B(g,0)$ as well.  It follows that $\mathbb D_a\subset B(g,0)$, and $g_n$ converges uniformly on compact subsets of $B(g,0)$.

\begin{figure}[h]
\input{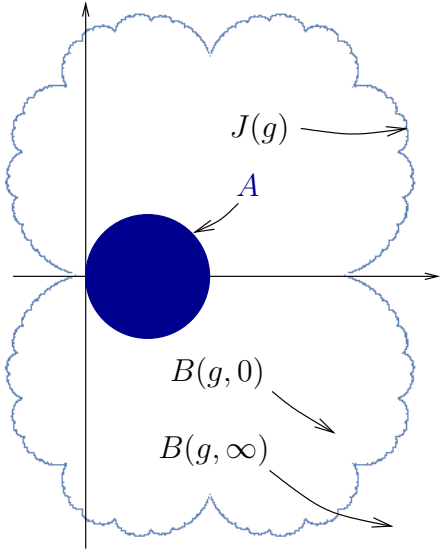}
\caption{$B(g,0)$ contains $A$.}\label{FIG}
\end{figure}

We complete the proof by adapting the proof of the Petal Theorem from \cite{beardon}, for which we need one more conjugacy.
Observe that $g$ is conjugate by $\sigma_0(z)= 1/z$ to $h\colon\mathbb C\rightarrow\mathbb C$ given by
\[h(z)=z+1+\frac{1}{z-1}.\]  
The fixed point $z=0$ for $g$ corresponds to the fixed point at $\infty$ for $h$.  It is easy to show that the half plane $\{z\in\mathbb C\colon\Re z>3\}$ is forward invariant by $h$, so it is contained in $B(h,\infty)$.  The image of this set by $\sigma_0$ is a disk of radius $1/6$ centered at $1/6$, contained in $B(g,0)$ (since the original set was contained in $B(h,\infty)$).  However, based on the picture of this basin in Figure \ref{FIG}, it appears we can construct a much larger forward invariant set.

Instead of starting with a half plane for $h$, we'll use the following parabolic region.  For each $t\geq0$, let
\[\Pi+t:=\{z=x+iy\colon y^2>12(3+t-x)\}.\]
It can be shown that $\Pi:=\Pi+0$ is forward invariant by $h$.  In particular, if $z\in\Pi$, we will show that $h(z)\in\Pi+1/2$.  Let $z=x+iy$, $h(z)=X+iY$, and $1/(z-1)=u+iv$, so $X=x+1+u$ and $Y=y+v$; then
\begin{eqnarray*}
Y^2-12(3+1/2-X)&=&(y+v)^2-12(3+1/2-(x+1+u))\\
&=&y^2-12(3+1/2-x)+v^2+2yv+12(1+u)\\
&>&v^2+2yv+12(1+u)\\
&\geq&2yv+12(1+u)\\
&\geq&12-2|yv|-12|u|\\
&>&0,
\end{eqnarray*}
where the last inequality can be derived from the fact that $z\in\Pi$ implies $|z|>3$.  Since $Y^2>12(3+1/2-X)$, we have $h(z)\in\Pi+1/2$.  We also have inductively that if $z\in\Pi$, then for all positive integers $n$,
\[h_n(z)\in\Pi+\frac{n}{2}.\]
Thus, for $h_n(x)+X_n+iY_n$,
\[|h_n(z)|^2=X_n^2+Y_n^2>X_n^2+12(2+n/2-X_n)=(X_n-6)^2+6n>n,\]
so $|h_n(z)|>\sqrt n$ and $h_n(z)\rightarrow\infty$ uniformly on $\Pi$.  The image of $\Pi$ by $\sigma_0$ is the cardioid
\[P:=\sigma_0(\Pi)=\{z=re^{i\theta}\colon6r<1+\cos\theta\},\]
and $g_n(z)\rightarrow0$ uniformly on $P$.  The set $P$ is the ``petal'' referred to in the Petal Theorem.

The preimage of $\mathbb D_a$ by $\sigma_0$ is the half plane
\[H_a:=\sigma_0(\mathbb D_a)=\{z=x+iy\colon 2\Re a+2\Im a>1\}.\]
Note that $\partial H_a$, the boundary of $H_a$, intersects the $x$-axis at $x=1/(2\Re a)$, and since $\Re a>0$, it is also never a horizontal line.  Moreover, $\partial\Pi$ is a horizontally oriented parabola, intersecting the $x$-axis at $x=3$.  See Figure \ref{FIG:PARAB}.    Thus, if $1/6<\Re a<1/4$, then $\partial H_a$ intersects the $x$-axis at $x<3$, so $\partial H_a$ must intersect $\partial\Pi$ at exactly two finite points.  If $0<\Re a\leq1/6$, then $\partial H_a$ intersects $\partial\Pi$ at exactly two finite points, one point (at which $\partial H_a$ is tangent to $\partial\Pi$), or zero points.  In the last two cases, we have $H_a\subset\Pi$. 

Returning to the coordinates centered at zero, this implies that
either $\mathbb D_a\subset P$ or $\partial\mathbb D_a$ intersects $\partial P$ at exactly two nonzero points.  
Recall that the set $\mathbb D_a$ is a disk of radius $|a|$ centered at $\bar a$, so both $\partial P$ and $\partial\mathbb D_a$ always contains the origin.  
\begin{figure}[h]
\input{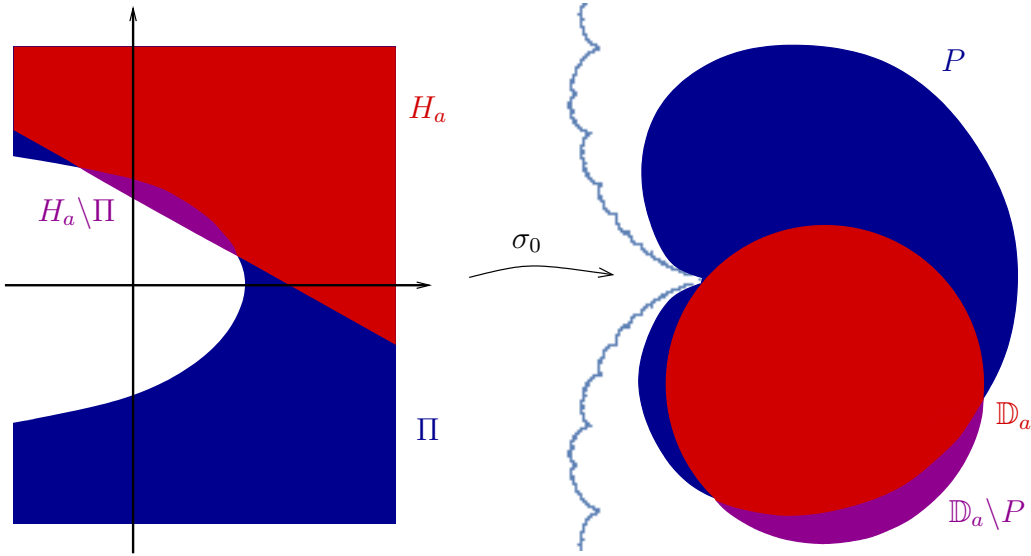}
\caption{On the left are $\Pi$, $H$, and $H_a\backslash\Pi$.  On the right are the images of these sets by $\sigma_0$: the cardioid, $P$, and $\mathbb D_a$ with boundaries intersecting at two nonzero points. }\label{FIG:PARAB}
\end{figure}
If $\partial\mathbb D_a$ intersects $\partial P$ at exactly two nonzero points, 
$\overline{\mathbb D}_a\backslash P$ is nonempty with a boundary comprised of the curve segment of $\partial P$ between the two nonzero intersections and the curve segment of $\partial\overline{\mathbb D}_a$ outside $P$ and between the two nonzero intersections.  Thus, $\overline{\mathbb D}_a\backslash P$ is closed and bounded, so it is compact.  It follows that $\overline{\mathbb D}_a\backslash P$ is always a strict, compact (sometimes trivially) subset of $B(g,0)$.

Since $g_n$ converges uniformly on $P$ and $\overline{\mathbb D}_a\backslash P\subset B(g,0)$, we have that $g_n$ converges uniformly on $\mathbb D_a$. Therefore, $\ph_n$ converges uniformly on \D. 

\end{proof}

Now we have the full picture: 

\begin{cor}\label{uci}
Suppose $\ph(z) = az^2+(1-2a)z+a$ and $|a-\frac{1}{4}| \leq \frac{1}{4}$. Then the iterates of \ph\ converge uniformly to $1$ on the entire open disk \D. 
\end{cor}

\section{\ph\ is essentially linear fractional}

We first provide the definition of an \textit{essentially linear fractional} map as given in \cite{bourdon2012spectra}. 

\begin{definition}[\cite{bourdon2012spectra}]\label{elf}
An analytic self-map \ph\ of \D\ is \textit{essentially linear fractional} if 

\begin{enumerate}
    \item \ph(\D) is contained in a proper subsdisk of \D\ internally tangent to the unit circle at some point $\eta \in \partial \D$;
    \item $\ph^{-1}(\{\eta\}) := \{ \gamma \in \partial \D : \eta \textrm{ belongs to the cluster set of } \eta \textrm{ of } \ph\ \textrm{ at } \gamma \}$ consists of one element, say $\zeta \in \partial \D$; and 
    \item $\ph'''$ extends continuously to $\D \cup \{ \zeta \}$. 
\end{enumerate}
\end{definition}

We will quickly verify that \textit{most} of the the maps $\ph(z) = az^2+(1-2a)z+a$ satisfy these conditions. (Again, we assume $a \neq 0$.) To do so, we need another result from \cite{bourdon2012spectra}.

\begin{proposition}[\cite{bourdon2012spectra}, Proposition 1.3]\label{derivatives}
Let \ph\ be an analytic self-map of \D\ that extends to be continuous on $\partial \D$. Suppose that $\ph \in C^2(1)$, that $\ph(1) = 1$, and that $|\ph(\zeta)| < 1$ for $\zeta \in \partial \D \backslash \{ 1 \}$. If 

$$\Re \left( \frac{1}{\ph'(1)} -1 + \frac{\ph''(1)}{\ph'(1)^2} \right) > 0 $$ then $\ph(\D)$ is contained in a proper subdisk of \D\ internally tangent to $\partial \D$ at $1$. 

\end{proposition}

Now to verify which our family of symbols \ph\ satisfy Definition \ref{elf}:

\begin{enumerate}
    \item For our symbols, the equation in Proposition \ref{derivatives} simplifies to $2a$, and we know $\Re a > 0$. However, Proposition \ref{derivatives} also requires that the rest of the unit circle is mapped into the \D. If $|a-\frac{1}{4}| < \frac{1}{4}$, Proposition \ref{derivatives} applies and $\ph(\D)$ is contained in a subdisk of \D\ internally tangent at $1$. However, we cannot use this proposition if $|a-\frac{1}{4}| = \frac{1}{4}$, since $-1$ maps to $4a-1$. If $a \neq \frac{1}{2}$, however, it is worth noting that $\ph \circ \ph$ satisfies the definition. 
    \item Since our function is analytic on the boundary, we are asking that $\ph^{-1}(\{ 1 \})$ contains only a single point from $\partial \D$. The points that map to $1$ are the zeroes of $\ph(z) - 1 $, $1$ and $1-\frac{1}{a}$, and the latter is outside of $\overline{\D}$ if $a \neq \frac{1}{2}$. If $a = \frac{1}{2}$, then \ph\ is \textit{not} essentially linear fractional since $-1$ also maps to $1$, and neither is any iterate of \ph, so we will handle that case separately.
    \item Since $\ph'''(z) \equiv 0$, this is trivial. 
\end{enumerate}

Unsurprisingly, just as with our work on UCI, we see that the defition of essentially linear fractional splits our work into the cases when $|a-\frac{1}{4}| < \frac{1}{4}$ and $|a-\frac{1}{4}| = \frac{1}{4}$. This continues in the following spectral results. 

\section{Spectrum of \W}

The uniform convergence of the iterates \ph\ on all of \D\ give the following spectral containment, found in \cite{cowen2016spectra}.

\begin{theorem}\cite[Corollary 10]{cowen2016spectra}\label{UCIspectra}.
Suppose $\ph:\D\rightarrow\D$ is
analytic with Denjoy-Wolff point $a$, $\ph_{n}\rightarrow a$ uniformly
in \D, and $\psi\in \Hi$ is continuous at $z=a$ with $\psi(a)\neq0$.
Then we have 
\[
\overline{\sigma_{p}(\psi(a)\C)}\subseteq\overline{\sigma_{ap}(T_{\psi}\C)}\subseteq\sigma(T_{\psi}\C)\subseteq\sigma(\psi(a)\C)
\]

In particular, if $\overline{\sigma_{p}(\C)}=\sigma(\C)$,
then $\sigma(T_{\psi}\C)=\sigma(\psi(a)\C)$.
\end{theorem}

Composition operators with parabolic symbols are notoriously more difficult when it comes to spectral problems, and this is no different. We have little information about $\sigma_{p}(\C)$; instead our goal is to use the fact that Theorem \ref{UCIspectra} gives us $\sigma(\W) \subseteq \sigma( \psi(a) \C) $. We now turn to two results from \cite{bourdon2012spectra} regarding essentially linear fractional maps. 

\begin{theorem}\cite[Theorem 3.3]{bourdon2012spectra}\label{bourdonCO}
Suppose that \ph\ is an essentially linear fractional self-map of \D\ fixing $1$. Suppose also that for $s =  \ph''(1)$, $\Re$ $ s > 0$. Then 
$$\sigma(\C) = \sigma_{e}(\C) = \{ e^{-s t} : t \geq 0 \} \cup \{ 0 \}. $$
\end{theorem}

\begin{theorem}\cite[Theorem 4.3]{bourdon2012spectra}\label{bourdonWCO}
Suppose \ph\ is essentially linear fractional with \ph(1) = 1, and $\psi \in \Hi$ is continuous at $1$. Then $\W \equiv \psi(1) \C$ modulo the compact operators.
\end{theorem}

Putting these facts together, we arrive at our main theorem.

\begin{theorem}\label{maintheorem1}
Suppose $\ph(z) = az^2+(1-2a)z+a $ maps \D\ into \D\, and $|a-\frac{1}{4}| < \frac{1}{4}$. Then for any $\psi \in \Hi$ continuous at $1$, we have 
$$\sigma(\W) = \sigma(\psi(1) \C) = \{ \psi(1) e^{-2a t} : t \geq 0 \} \cup \{ 0 \}. $$
\end{theorem}
\begin{proof}
By Theorem \ref{UCIspectra}, we have $\sigma(\W) \subseteq \sigma(\psi(1) \C) $. By Theorem \ref{bourdonCO}, we have $\sigma(\C) = \sigma_{e}(\C)$. By Theorem \ref{bourdonWCO}, we have $\sigma_{e}(\W) = \sigma_{e}(\psi(1)\C)$. Putting these together, we have
$$\sigma(\psi(1) \C) = \sigma_{e}(\psi(1)\C) = \sigma_{e}(\W) \subseteq \sigma(\W) \subseteq \sigma(\psi(1)\C),$$

\noindent and since the first and last sets in the containment are equal, we have $\sigma(\psi(1)\C) = \sigma(\W)$. By Theorem \ref{bourdonCO}, noting $\ph''(1) = 2a$, we have $$\sigma(\W) = \{ \psi(1) e^{-2a t} : t \geq 0 \} \cup \{ 0 \}. $$
\end{proof}

We now investigate the situation when $|a-\frac{1}{4}|=\frac{1}{4}$. While we have shown that \ph\ still converges under iteration to $1$ uniformly on all of \D, it is \textit{not} essentially linear fractional since $|\ph(-1)| = 4a-1| = 1$. Here, we must actually divide our special case yet again: if $a\neq \frac{1}{2}$, then $\ph_2 = \ph \circ \ph$ is essentially linear fractional (since $4a-1$ is then mapped into \D) and of course $\ph_2$ is uniformly convergent under iteration on all of \D. Therefore, we as a corollary to Theorem \ref{maintheorem1}, we get the following:

\begin{cor}\label{maintheorem2}
Suppose $\ph(z) = az^2+(1-2a)z+a $ maps \D\ into \D\, and $|a-\frac{1}{4}| = \frac{1}{4}$, $a \neq 0, \frac{1}{2}$. Then for any $\psi \in \Hi$ continuous at $1$, we have 

$$\sigma(W_{\psi,\ph_2}) = \sigma(\psi(1) C_{\ph_{2}}) = \{ \psi(1) e^{-4a t} : t \geq 0 \} \cup \{ 0 \} $$
where $\ph_2 = \ph \circ \ph$.
\end{cor}

\begin{proof}
The proof follows exactly as before, except that $\ph_2''(1) = 4a$. 
\end{proof}

The result of Corollary \ref{maintheorem2} suggests that it is most likely true that $\sigma(\W)$ is the same as shown in \ref{maintheorem1} when $|a-\frac{1}{4}|=\frac{1}{4}$, but we do not have a proof. 

However, even then, $a = \frac{1}{2}$ proves itself to be an entirely distinct case. Here, instead of functional behavior, we now rely on a linear algebra trick also used in \cite{cowen2016spectra}.

\begin{lemma}\cite[Lemma 3]{cowen2016spectra}\label{ab}
If $A$ and $B$ are bounded linear operators on a Hilbet space $ \mathcal{H}$, then $\sigma(AB) \cup \{ 0 \} = \sigma(BA) \cup \{ 0 \}$.
\end{lemma}

Using this, we can now finish the story with $a = \frac{1}{2}$, which actually varies just slightly from the result in \ref{maintheorem1}. 

\begin{theorem}\label{specialcase}
Suppose $\ph(z) = \frac{1}{2}z^2+\frac{1}{2} $, an analytic self-map of \D. Then for any $\psi \in \Hi$ continuous at $1$, we have 

$$\sigma(\W) = \{ \psi(1) e^{- t/2} : t \geq 0 \} \cup \{ 0 \}. $$
\end{theorem}
\begin{proof}
If $f(z) = \frac{1}{2}z+\frac{1}{2}$, then we have $\C = C_{z^2}C_f$. By Lemma \ref{ab}, we have $\sigma(\C) \cup \{ 0 \} = \sigma(C_{z^2}C_f) \cup \{ 0 \} = \sigma(C_f C_{z^2}) \cup \{ 0 \} = \sigma(C_{f^2}) \cup \{ 0\} $. Note that $f^2(z) = \frac{1}{4}z^2+\frac{1}{2}z+\frac{1}{4}$, which falls under Theorem \ref{bourdonCO} with $s = \frac{1}{2}$. Since we also know \C\ is not invertible, we have $\sigma(\C) = \{ \psi(1) e^{-t/2} : t \geq 0 \} \cup \{ 0 \}. $

Likewise, let $\psi \in \Hi$ be continuous at $1$ and consider $\T \C = \T  C_{z^2}C_f$. Again by Lemma \ref{ab}, we have $\sigma(\T  C_{z^2}C_f) \cup \{ 0 \} = \sigma( C_f \T  C_{z^2}) \cup \{ 0 \} = \sigma(T_{\psi \circ f} C_{f^2}) \cup \{ 0 \} $. Since $f$ maps \D\ into \D\ analytically and $f$ is continuous at $1$, we still have that $\psi \circ f \in \Hi$ and and $\psi \circ f$ is continuous at $1$ (and $\psi \circ f (1) = \psi(1)$). Again, we also know that $\T \C$ is not invertible. Then, by Theorem \ref{maintheorem1}  we have

$$\sigma(\Ws) = \sigma(\psi(1) \C) = \{ \psi(1) e^{-t/2} : t \geq 0 \} \cup \{ 0 \}. $$
\end{proof}

While our guess is that Theorem \ref{maintheorem1} holds for $|a-\frac{1}{4}|=\frac{1}{4}$ when $a$ is complex, the exponent in our result for $a = \frac{1}{2}$ does not match up with Theorem \ref{maintheorem1}, presumably because it is more distinct in its failure to be essentially linear fractional. However, the final result is the same in practicality: for $0 < a \leq \frac{1}{2}$, $\sigma(\W)$ is the closed line segment connecting $\psi(1)$ to the origin.

\section{Implications and Further Questions}

There are two important concepts in this paper that could be pursued further. The first is extending our methods for showing quadratics exhibit uniformly convergent iteration to higher-degree polynomials. Certainly Corollary \ref{uci} implies more than it says; e.g. $\ph \circ \ph$ is a quartic self-map of \D\ that also converges uniformly on all of \D. For any polynomial map fixing $1$, $\ph_{n}(z)-1$ will be a factor of $\ph_{n+1}(z)-1$, suggesting that our recursive approach used in Theorem \ref{uci1} could be generalized, but it will require deeper geometric intuition than we use here for quadratic maps. 
    
The second important concept is the intersection of self-maps of \D\ that exhibit uniformly convergent iteration on all of \D, and essentially linear fractional maps. Certainly they do not perfectly align; we have already seen that $\frac{1}{2}z^2+\frac{1}{2}$ converges under iteration to $1$ uniformly on all of \D, but is not essentially linear fractional. Likewise, a linear fractional map with both an interior and a boundary fixed point (e.g. $\ph(z) = \frac{z}{2-z}$) cannot converge uniformly on all of \D; it must have only one fixed point in $\overline{\D}$ \cite[Theorem 3]{cowen2016spectra}. However, the intersection of the two concepts is non-trivial, and leads to the following conjecture.

\begin{conj}
Suppose \ph\ is an essentially linear fractional self-map of \D\ with exactly one fixed point $w$ in $\overline{\D}$. Then \ph\ converges under iteration to $w$ uniformly on all of \D. 
\end{conj}

Were this conjecture true, it would immediately provide a full description of $\sigma(\W)$ for a rather broad class of symbols, by the same arguments made in this paper. Thus we end with the following list of questions:

\begin{enumerate}
    \item What is $\sigma(\W)$ if $\ph(z)=az^2+(1-2a)z+a$ and $|a-\frac{1}{4}| = \frac{1}{4}, a \neq 0, \frac{1}{2}$ ? 
    \item When do essentially linear fractional maps converge uniformly to their Denjoy-Wolff point on all of \D?
    \item Which polynomial self-maps of \D\ converge uniformly on all of \D\ to the Denjoy-Wolff point?
    \item If $\psi \in \Hi$ is continuous at the Denjoy-Wolff point $w$ of \ph, and \ph\ is not an automorphism, then when, if ever, is $\sigma(\W) \neq \sigma(\psi(w)\C)?$
\end{enumerate}

\section*{Acknowledgements}

The authors would like to thank NSF for funding this research, and the directors of CURM (Kathryn Leonard, Maria Mercedes Franco) for their counsel. We would also like to thank Michal Misiurewicz, Paul Bourdon, and Carl Cowen for their help with certain difficulties in the proofs. 
\footnotesize

{\footnotesize{
    \bibliography{citations}
    \bibliographystyle{unsrt}
    }}

\end{document}